# TUSNÁDY'S INEQUALITY REVISITED

By Andrew Carter and David Pollard

*University of California, Santa Barbara and Yale University*

Tusnády's inequality is the key ingredient in the KMT/Hungarian coupling of the empirical distribution function with a Brownian bridge. We present an elementary proof of a result that sharpens the Tusnády inequality, modulo constants. Our method uses the beta integral representation of Binomial tails, simple Taylor expansion and some novel bounds for the ratios of normal tail probabilities.

**1. Introduction.** In one of the most important probability papers of the last forty years, Komlós, Major and Tusnády (1975) sketched a proof for a very tight coupling of the standardized empirical distribution function with a Brownian bridge, a result now often referred to as the KMT, or Hungarian, construction. Their coupling greatly simplifies the derivation of many classical statistical results—see Shorack and Wellner [(1986), Chapter 12 et seq.], for example.

The construction has taken on added significance for statistics with its use by Nussbaum (1996) in establishing asymptotic equivalence of density estimation and white noise models. Brown, Carter, Low and Zhang (2004) have somewhat simplified and expanded Nussbaum's argument using our Theorem 2, via inequality (5).

At the heart of the KMT method [with refinements as in the exposition by Csörgő and Révész (1981), Section 4.4] lies the quantile coupling of the $\text{Bin}(n, 1/2)$ and $N(n/2, n/4)$ distributions, which may be defined as follows. Let $Y$ be a random variable distributed $N(n/2, n/4)$. Find the cutpoints $-\infty = \beta_0 < \beta_1 < \cdots < \beta_n < \beta_{n+1} = \infty$ for which

$$\mathbb{P}\{\text{Bin}(n, 1/2) \geq k\} = \mathbb{P}\{Y > \beta_k\} \qquad \text{for } k = 0, 1, \ldots, n.$$

When $\beta_k < Y \le \beta_{k+1}$, let $X$ take the value $k$. Then $X$ has a $\text{Bin}(n, 1/2)$ distribution.

It is often more convenient to work with the tails of the standard normal $\bar{\Phi}(z) = \mathbb{P}\{N(0,1) > z\}$ and the standardized cutpoint $z_k = 2(\beta_k - n/2)/\sqrt{n}$, thereby replacing $\mathbb{P}\{Y > \beta_k\}$ by $\bar{\Phi}(z_k)$.

Symmetry considerations show that $\beta_{n-k+1} = n - \beta_k$, so that it suffices to consider only half the range for $k$. More precisely, when $n$ is even, say $n = 2m$, the interval $(\beta_m, \beta_{m+1})$ is symmetric about $n/2$, so we have only to consider $k \ge m + 1 = (n+2)/2$. When $n$ is odd, say $n = 2m+1$, the interval $(\beta_m, \beta_{m+2})$ is symmetric about $n/2 = \beta_{m+1}$, so we have only to consider $k \ge m + 2 = (n+3)/2$.

The usual normal approximation with continuity correction suggests that $\beta_k \approx k - 1/2$, which, if true, would bound $|X - Y|$ by a constant that does not change with $n$. Of course, such an approximation for all $k$ is too good to be true, but results almost as good have been established. The most elegant version appeared in the unpublished dissertation (in Hungarian) of Tusnády (1977), whose key inequality may be expressed as the assertion

$$(1) \qquad k - 1 \le \beta_k \le \frac{3n}{2} - \sqrt{2n(n-k)} \qquad \text{for } n/2 \le k \le n.$$

As explained by Csörgő and Révész [(1981), Section 4.4], Tusnády's inequality implies that $|X - n/2| \le |Y - n/2| + 1$ and $|X - Y| \le 1 + Z^2/8$, where $Z$ denotes the standardized variable $(2Y - n)/\sqrt{n}$. They also noted that Tusnády's proof of inequality (1) was "elementary," but "not at all simple." Bretagnolle and Massart [(1989), Appendix] published another proof of Tusnády's inequality—an exquisitely delicate exercise in elementary calculus and careful handling of Stirling's formula to approximate individual Binomial probabilities. With no criticism intended, we note that their proof is quite difficult. More recently, Dudley [(2000), Chapter 1] and Massart (2002) have reworked and refined the Bretagnolle/Massart calculations. Clearly, there is a continuing perceived need for an accessible treatment of the coupling result that underlies the KMT construction.

With this paper we offer another approach, which actually leads to an improvement (modulo constants) of the Tusnády inequality. In fact, the Tusnády upper bound greatly overestimates $\beta_k$ for moderate to large $k$. (See below.) Our method differs from that of Bretagnolle and Massart, in that we work directly with the whole tail probability. Our method is closer to that of Peizer and Pratt (1968), who suggested a Cornish–Fisher expansion of the Binomial percentiles—but, as noted by Pratt [(1968), Sections 5 and 8], a rigorous proof by this method is difficult. To avoid the difficulty, Molenaar [(1970), Section III.2] made a more direct calculation starting from the representation of the Binomial tail as a beta integral,

$$(2) \qquad \mathbb{P}\{\text{Bin}(n, 1/2) \ge k\} = \frac{n!}{(k-1)!(n-k)!} \int_0^{1/2} t^{k-1}(1-t)^{n-k}\, dt.$$

He indicated that his expansion would be valid provided $|k - n/2| = O(\sqrt{n})$. Pratt seemed to be claiming validity for his expansion for the range $|k - n/2| = o(n)$, but we believe extra work is needed for $|k - n/2|$ large.

We should point out that Peizer, Pratt and Molenaar were actually concerned with normal approximations to distributions more general than the $\mathrm{Bin}(n, 1/2)$ case needed for the KMT construction. We have specialized their results to this case.

Our method also starts from the integral representation (2), to derive an approximation via Laplace's method for integrals [de Bruijn (1981), Section 4.3] using only Taylor's theorem and Stirling's formula [Feller (1968), Section II.9]

$$
(3) \quad n! = \sqrt{2\pi} \exp((n + \tfrac{1}{2})\log n - n + \lambda_n)
$$
$$
\text{with } (12n+1)^{-1} \leq \lambda_n \leq (12n)^{-1}.
$$

In fact [Komlós, Major and Tusnády (1975), page 130], the KMT construction only needs a result like the Tusnády inequality for values of $k$ in a range where $|2k - n| \leq \varepsilon_0 n$ for some fixed $\varepsilon_0 < 1$. For that range, a suitable bound can be derived from classical large deviation approximations for Binomial tails. For example, in an expanded version of the argument sketched in the 1975 paper, Major (2000) used the large deviation approximation

$$\mathbb{P}\{X \geq k\} = \bar{\Phi}(\varepsilon\sqrt{n})\exp(A_n(\varepsilon)) \quad \text{where } \varepsilon = (2k-n)/n,$$

with

$$|A_n(\varepsilon)| = O(n\varepsilon^3 + n^{-1/2}) \quad \text{uniformly in } 0 \leq \varepsilon \leq \varepsilon_0 < 1.$$

Mason (2001) derived the KMT coupling from an analogous approximation with

$$A_n(\varepsilon) = n\varepsilon^3 \lambda(\varepsilon) + O(\varepsilon + n^{-1/2}) \quad \text{uniformly in } 0 \leq \varepsilon \leq \varepsilon_0 < 1,$$

where $\lambda(\cdot)$ is a power series whose coefficients depend on the cumulants of the Binomial distribution. Such an approximation follows from a minor variation on the general method explained by Petrov [(1975), Section 8.2]. Symmetry of the $\mathrm{Bin}(n, 1/2)$ makes the third cumulant zero; the power series $\varepsilon^3 \lambda(\varepsilon)$ starts with a multiple of $\varepsilon^4$.

Our method gives a sharper approximation to the $\mathrm{Bin}(n, 1/2)$ tails over the range $n/2 < k \leq n - 1$ (which, by symmetry, actually covers the range $0 < k < n$). Only at the extreme, $k = n$, does the calculation fail.

THEOREM 1. *Let $X$ have a $\mathrm{Bin}(n, 1/2)$ distribution, with $n \geq 28$. Define*

$$\gamma(\varepsilon) = \frac{(1+\varepsilon)\log(1+\varepsilon) + (1-\varepsilon)\log(1-\varepsilon) - \varepsilon^2}{2\varepsilon^4} = \sum_{r=0}^{\infty} \varepsilon^{2r}/(2r+3)(2r+4),$$

an increasing function with $\gamma(0) = 1/12$ and $\gamma(1) = -1/2 + \log 2 \approx 0.1931$. Define $\varepsilon = (2K - N)/N$, where $K = k - 1$ and $N = n - 1$. Define $\lambda_n$ as in (3). Then there is a constant $C$ such that

$$\mathbb{P}\{X \geq k\} = \bar{\Phi}(\varepsilon\sqrt{N})\exp(A_n(\varepsilon)),$$

where

$$A_n(\varepsilon) = -N\varepsilon^4\gamma(\varepsilon) - \tfrac{1}{2}\log(1-\varepsilon^2) - \lambda_{n-k} + r_k \quad and \quad -C\log N \leq Nr_k \leq C$$

for all $\varepsilon$ corresponding to the range $n/2 < k \leq n - 1$.

Notice that the $\lambda_{n-k}$ can be absorbed into the error terms, and that $\log(1 - \varepsilon^2)$ is small compared with $N\varepsilon^4 + O(n^{-1})$, when $\varepsilon \leq \varepsilon_0 < 1$.

A very precise approximation for the cutpoints $\beta_k$ follows from Theorem 1 inequalities (see Section 3) for the tails of the normal distribution.

THEOREM 2. *Let $z_k = 2(\beta_k - n/2)/\sqrt{n}$ and $\varepsilon = (2K - N)/N$. Let $S(\varepsilon) = \sqrt{1 + 2\varepsilon^2\gamma(\varepsilon)}$ for $\gamma(\varepsilon)$, as in Theorem 1. Then, for some constant $C'$ and $n \geq 28$,*

$$z_k = \varepsilon\sqrt{N}\,S(\varepsilon) + \frac{\log(1-\varepsilon^2) + 2\lambda_{n-k}}{2\varepsilon\sqrt{N}\,S(\varepsilon)} + \theta_k$$

*with $-C'(\varepsilon\sqrt{N} + 1) \leq N\theta_k \leq C'(\varepsilon\sqrt{N} + \log N)$ for all $\varepsilon$ corresponding to the range $n/2 < k \leq n - 1$.*

For example, the theorem implies $\beta_k - k + 1/2 = o(1)$ uniformly over a range where $|k - n/2| = o(n^{2/3})$. Also, when $\varepsilon \leq \varepsilon_0 < 1/2$, the log term can be absorbed into the $O(\varepsilon/\sqrt{n})$ errors. Even when $k$ gets close to $n - 1$, the log term contributes only an $O(n^{-1/2}\log n)$ to the approximation. More precisely, if $k = n - B$ for a fixed $B \geq 1$, our approximation simplifies to

$$(4) \quad \beta_{n-B} = \frac{1+c}{2}n - \frac{1+2B}{4c}\log n + O(1) \quad \text{where } c = S(1) \approx 1.177,$$

which agrees up to $O(1)$ terms with the result obtained by direct calculation from

$$\mathbb{P}\{X \geq n - B\} = \left(\binom{n}{0} + \cdots + \binom{n}{B}\right)2^{-n} = \frac{n^B}{B!2^n}(1 + o(1))$$

and the well-known approximation for normal percentiles,

$$\bar{\Phi}^{-1}(p) = y - \frac{\log y}{y} + O(1/y) \quad \text{as } p \to 0, \text{ where } y = \sqrt{2\log(1/p)}.$$

By contrast, the upper bound for $\beta_{n-B}$ from (1) is about $0.088n$ too large.

It is also an easy consequence of Theorem 2 that there exist positive constants $C_i$ for which

$$(5) \quad -\frac{C_1}{\sqrt{n}} + C_2 \frac{|k-n/2|^3}{n^2} \leq \beta_k - k + \frac{1}{2} \leq \frac{C_3 \log n}{\sqrt{n}} + C_4 \frac{|k-n/2|^3}{n^2}$$

for $n/2 \leq k \leq n$ and all $n$. For the quantile coupling between an $X$ distributed $\mathrm{Bin}(n, 1/2)$ and a $Y = n/2 + \sqrt{n}Z/2$ distributed $N(n/2, n/4)$, it follows that there is a positive constant $C$ for which

$$\left|X - \frac{n}{2}\right| \leq C + \left|Y - \frac{n}{2}\right| \quad \text{and} \quad |X - Y| \leq C + \frac{C}{n^2}\left|X - \frac{n}{2}\right|^3.$$

Using the fact that $|X - n/2| \leq n/2$, we could also write the upper bound for $|X - Y|$ as a constant multiple of $1 + Z^2(1 \wedge |Z|/\sqrt{n})$, which improves on Tusnády's $1 + Z^2/8$, modulo multiplicative constants. (We have made no attempt to find the best constants, even though, in principle, explicit values could be found by our method.)

**2. Outline of our method.** As in Theorem 1, write $\varepsilon = (2K - N)/N$, where $K = k - 1$ and $N = n - 1$. Then $K/N = (1 + \varepsilon)/2$ and the range $n/2 < k < n$ corresponds to

$$(6) \quad 1 - \frac{2}{N} \geq \varepsilon = \frac{2K}{N} - 1 \geq \begin{cases} N^{-1} & \text{when } n \text{ is even,} \\ 2N^{-1} & \text{when } n \text{ is odd.} \end{cases}$$

Define $2H(t) = (1+\varepsilon)\log t + (1-\varepsilon)\log(1-t)$ for $0 < t < 1$. Representation (2) can then be rewritten as

$$\mathbb{P}\{X \geq k\} = \frac{nN!}{K!(N-K)!} \int_0^{1/2} \exp(K \log t + (N-K)\log(1-t))\, dt$$

$$= \frac{nN!}{K!(N-K)!} \int_0^{1/2} e^{NH(t)}\, dt.$$

By Stirling's formula (3),

$$\frac{N!}{K!(N-K)!} = \frac{1}{N}\sqrt{\frac{4N}{2\pi(1-\varepsilon^2)}} \exp(\Lambda - NH(K/N))$$

where $\Lambda := \lambda_N - \lambda_K - \lambda_{N-K}$.

Thus, the beta integral equals

$$\frac{n}{N}\exp\left(\Lambda - \frac{1}{2}\log(1-\varepsilon^2) - NH(K/N)\right)\sqrt{\frac{4N}{2\pi}} \int_0^{1/2} e^{NH(t)}\, dt.$$

The function $H(\cdot)$ is concave on $(0, 1)$. It achieves its global maximum at $K/N$, which lies outside the range of integration. On the interval $(0, 1/2]$ the maximimum is achieved at $1/2$. On the range of integration, $H(t) - H(K/N)$ is never greater than

$$H(1/2) - H(K/N) = -\tfrac{1}{2}(1+\varepsilon)\log(1+\varepsilon) - \tfrac{1}{2}(1-\varepsilon)\log(1-\varepsilon)$$
$$= -\tfrac{1}{2}\varepsilon^2 - \varepsilon^4\gamma(\varepsilon).$$

The concave function $h(s) := H((1-s)/2) - H(1/2)$ achieves its maximum value of zero at $s=0$ and

$$\mathbb{P}\{X \geq k\} = e^\Delta \sqrt{\frac{N}{2\pi}} \int_0^1 e^{Nh(s) - N\varepsilon^2/2} \, ds, \tag{7}$$

where $\Delta = \log(1 + N^{-1}) + \Lambda - \tfrac{1}{2}\log(1-\varepsilon^2) - N\varepsilon^4\gamma(\varepsilon)$.

The $\Delta$ contributes $O(1/n) - \lambda_{n-k} - \tfrac{1}{2}\log(1-\varepsilon^2) - N\varepsilon^4\gamma(\varepsilon)$ to the $A_n(\varepsilon)$ from Theorem 1. Taylor's expansion of $h(s)$ about $s=0$ and concavity of $h(\cdot)$ show that the exponent $Nh(s)$ drops off rapidly as $s$ moves away from zero. Indeed,

$$\begin{aligned} h(s) &= -\varepsilon s - \tfrac{1}{2}s^2 + \tfrac{1}{6}s^3 h'''(s^*) && \text{with } 0 < s^* < s \\ &\approx \tfrac{1}{2}\varepsilon^2 - \tfrac{1}{2}(s+\varepsilon)^2 && \text{for } s \text{ near zero.} \end{aligned} \tag{8}$$

See Section 4 for the more precise statement of the approximation.

Most of the contribution to the integral (7) comes from $s$ in a small neighborhood of 0. Ignoring tail contributions to the integral, we will then have

$$\mathbb{P}\{X \geq k\} \approx e^\Delta \sqrt{\frac{N}{2\pi}} \int_0^\infty \exp\left(-\tfrac{1}{2}N(s+\varepsilon)^2\right) ds = e^\Delta \bar{\Phi}(\varepsilon\sqrt{N}), \tag{9}$$

as asserted by Theorem 1.

To derive Theorem 2 we perturb the argument $\varepsilon\sqrt{N}$ slightly to absorb the factor $\exp(A_n(\varepsilon))$. We seek a $y$ for which

$$\bar{\Phi}(\varepsilon\sqrt{N} + y) \approx \exp(A_n(\varepsilon))\bar{\Phi}(\varepsilon\sqrt{N}) = \bar{\Phi}(z_k).$$

That is, we need

$$\bar{\Phi}(\varepsilon\sqrt{N} + y)/\bar{\Phi}(\varepsilon\sqrt{N}) \approx \exp(-N\varepsilon^4\gamma(\varepsilon) - \tfrac{1}{2}\log(1-\varepsilon^2)).$$

As shown in the next section, the ratio of normal tail probabilities $\bar{\Phi}(x+y)/\bar{\Phi}(x)$ behaves like $\exp(-xy - y^2/2)$, at least when $x$ is large. Ignore the logarithmic term for the moment. Then the heuristic suggests that we choose $y$ to make $\varepsilon\sqrt{N}y + y^2/2 \approx N\varepsilon^4\gamma(\varepsilon)$, that is,

$$y \approx -\varepsilon\sqrt{N} + \sqrt{N\varepsilon^2 + 2N\varepsilon^4\gamma(\varepsilon)}$$

and, hence,

$$z_k \approx \varepsilon\sqrt{N} + y \approx \varepsilon\sqrt{N}\sqrt{1 + 2\varepsilon^2\gamma(\varepsilon)}.$$

For the rigorous proof of Theorem 2 we need to replace these heuristic approximations by inequalities giving upper and lower bounds for $\bar{\Phi}(z_k)$, then invoke the inequalities for normal tails derived in the next section.

**3. Tails of the normal distributions.** The classical tail bounds for the normal distribution [cf. Feller (1968), Section VII.1 and Problem 7.1] show that $\bar{\Phi}(x)$ behaves roughly like the density $\phi(x)$:

(10)
$$\left(\frac{1}{x} - \frac{1}{x^3}\right)\phi(x) < \bar{\Phi}(x) < \frac{1}{x}\phi(x) \qquad \text{for } x > 0.$$
$$\bar{\Phi}(x) < \tfrac{1}{2}\exp(-x^2/2)$$

The first upper bound is good for large $x$, the second for $x \approx 0$. For the proofs of both Theorem 1 (in Section 4) and Theorem 2 (in Section 5), we will need to bound the ratio $\bar{\Phi}(x+y)/\bar{\Phi}(x)$. It is possible to derive suitable bounds directly from (10), but we have found it easier to work with inequalities that interpolate smoothly between the different cases in (10). We express our results in logarithmic form, using the function $\Psi(x) := -\log\bar{\Phi}(x)$ and its derivative

$$\rho(x) = \frac{d}{dx}\Psi(x) = \phi(x)/\bar{\Phi}(x).$$

To a first approximation, the positive function $\rho(x)$ increases like $x$. By inequality (10), the error of approximation, $r(x) := \rho(x) - x$, is positive for $x > 0$ and, for $x > 1$,

$$r(x) < \frac{x}{x^2 - 1} = O(1/x) \qquad \text{as } x \to \infty.$$

In fact, as shown by the proof of the next lemma, $\rho(\cdot)$ is increasing and $r(\cdot)$ is decreasing and positive on the whole real line.

LEMMA 1. *The function $\rho(\cdot)$ is increasing and the function $r(\cdot)$ is decreasing, with $r(\infty) = \rho(-\infty) = 0$ and $r(0) = \rho(0) = 2/\sqrt{2\pi} \approx 0.7979$. For all $x \in \mathbb{R}$ and $\delta \geq 0$, the increments of the function $\Psi(x) := -\log\bar{\Phi}(x)$ satisfy the following inequalities:*

(i) $\delta\rho(x) \leq \Psi(x+\delta) - \Psi(x) \leq \delta\rho(x+\delta),$
(ii) $\delta r(x+\delta) \leq \Psi(x+\delta) - \Psi(x) - \tfrac{1}{2}(x+\delta)^2 + \tfrac{1}{2}x^2 \leq \delta r(x),$
(iii) $x\delta + \tfrac{1}{2}\delta^2 \leq \Psi(x+\delta) - \Psi(x) \leq \rho(x)\delta + \tfrac{1}{2}\delta^2.$

PROOF. Let $Z$ be $N(0,1)$ distributed. Define $M(x) = \mathbb{P}e^{-x|Z|}$, a decreasing function of $x$ with $\log M(x)$ strictly convex. Notice that

$$1/\rho(x) = \sqrt{2\pi} \exp(x^2/2) \int_0^\infty \phi(z+x)\,dz$$

$$= \int_0^\infty \exp(-xz - z^2/2)\,dz = \sqrt{\frac{\pi}{2}} M(x).$$

Thus, $-\log M(x) - \log\sqrt{\pi/2} = \log \rho(x) = \Psi(x) - x^2/2 - \log\sqrt{2\pi}$ is a concave, increasing function of $x$ with derivative $\rho(x) - x = r(x)$. It follows that $r(\cdot)$ is a decreasing function, because

$$r'(x) = -\frac{d^2}{dx^2}\log M(x) < 0 \qquad \text{by convexity of } \log M(x).$$

Inequality (i) follows from the equality

$$\Psi(x+\delta) - \Psi(x) = \delta\Psi'(y^*) = \delta\rho(y^*) \qquad \text{for some } x < y^* < x+\delta,$$

together with the fact that $\rho(\cdot)$ is an increasing function. Similarly, the fact that

$$\frac{d}{dy}\left(\Psi(y) - \frac{1}{2}y^2\right) = \rho(y) - y = r(y) \qquad \text{which is a decreasing function}$$

gives inequality (ii). Inequality (iii) follows from (ii) because $\delta r(x+\delta) \geq 0$ and $x\delta + r(x)\delta = \rho(x)\delta$. □

Reexpressed in terms of the tail function $\bar\Phi$, the three inequalities from the lemma become:

(i) $\exp(-\delta\rho(x)) \geq \bar\Phi(x+\delta)/\bar\Phi(x) \geq \exp(-\delta\rho(x+\delta))$,
(ii) $\exp(-\delta r(x+\delta)) \geq \exp(x\delta + \delta^2/2)\bar\Phi(x+\delta)/\bar\Phi(x) \geq \exp(-\delta r(x))$,
(iii) $\exp(-x\delta - \delta^2/2) \geq \bar\Phi(x+\delta)/\bar\Phi(x) \geq \exp(-\rho(x)\delta - \delta^2/2)$.

Less formally,

$$\mathbb{P}\{Z \leq x+\delta \mid Z \leq x\} = 1 - \bar\Phi(x+\delta)/\bar\Phi(x) \approx \delta\rho(x) \qquad \text{for small } \delta,$$

which corresponds to the fact that $\rho$ is the hazard rate for the $N(0,1)$ distribution.

**4. Details of the proof for Theorem 1.** To make the proof rigorous, we need to replace the approximation in the Taylor expansion (8) by upper and lower bounds involving the third derivative

$$h'''(s) = \frac{1-\varepsilon}{(1+s)^3} - \frac{1+\varepsilon}{(1-s)^3} = -\frac{6s + 2s^2 + \varepsilon(2+6s^2)}{(1-s^2)^3}.$$

The derivative of this function is negative for all $s$. Thus,

$$h'''(s) \leq h'''(0) = -2\varepsilon \qquad \text{for } 0 < s < 1$$

and

$$h(s) \leq \tfrac{1}{2}\varepsilon^2 - \tfrac{1}{2}(s+\varepsilon)^2 \qquad \text{for } 0 < s < 1.$$

The right-hand side of the approximation (9) is actually an upper bound, because the integrand is nonnegative on $(1, \infty)$. That is,

$$\mathbb{P}\{X \geq k\} \leq e^\Delta \bar{\Phi}(\varepsilon\sqrt{N}),$$

which gives the upper bound for $A_n(\varepsilon)$ stated in the theorem.

For the lower bound, for some small positive $\eta$ discard the contribution to the integral in (7) from the range $(\eta, 1)$, and bound $h'''$ from below by $h'''(\eta)$ on the range $(0, \eta)$, then integrate to get

$$\mathbb{P}\{X \geq k\} \geq e^\Delta \sqrt{\frac{N}{2\pi}} \int_0^\eta \exp\left(-\frac{1}{2}N(s+\varepsilon)^2 + \frac{1}{6}N\eta s^2 h'''(\eta)\right) ds$$

$$= e^\Delta \sqrt{\frac{N}{2\pi}} \int_0^\eta \exp\left(-\frac{1}{2}N\kappa^2(s+\varepsilon/\kappa^2)^2 + \frac{1}{2}N\varepsilon^2/\kappa^2 - \frac{1}{2}N\varepsilon^2\right) ds$$

$$= \frac{e^\Delta}{\kappa} \exp\left(\frac{1}{2}N\varepsilon^2/\kappa^2 - \frac{1}{2}N\varepsilon^2\right)(\bar{\Phi}(\varepsilon\sqrt{N}/\kappa) - \bar{\Phi}(\varepsilon\sqrt{N}/\kappa + \kappa\eta\sqrt{N})),$$

where

$$\kappa^2 = 1 - \tfrac{1}{3}\eta h'''(\eta) \leq 1 + 6\eta(\eta + \varepsilon) \qquad \text{if } \eta \leq \tfrac{1}{2}.$$

From Lemma 1, parts (iii) and (ii),

$$\bar{\Phi}(\varepsilon\sqrt{N}/\kappa + \kappa\eta\sqrt{N}) \leq \bar{\Phi}(\varepsilon\sqrt{N}/\kappa)\exp(-N\varepsilon\eta - \tfrac{1}{2}N\kappa^2\eta^2)$$

and

$$\exp(\tfrac{1}{2}N\varepsilon^2)\bar{\Phi}(\varepsilon\sqrt{N}) \leq \exp(\tfrac{1}{2}N\varepsilon^2/\kappa^2)\bar{\Phi}(\varepsilon\sqrt{N}/\kappa).$$

Thus

(11) $\quad \mathbb{P}\{X \geq k\} \geq \exp(\Delta - \log\kappa)\bar{\Phi}(\varepsilon\sqrt{N})[1 - \exp(-N\varepsilon\eta - \tfrac{1}{2}N\kappa^2\eta^2)].$

We need $\log \kappa = O(\ell_N)$, where $\ell_N = N^{-1}\log N$, for otherwise the asserted inequality $-C\log N \leq r_k$ would be violated. As $\log \kappa \leq 6(\eta^2 + \eta\varepsilon)$, this requirement suggests that we take $\eta$ as a solution to the equation $\tfrac{1}{2}\eta^2 + \eta\varepsilon = \ell_N$, that is, $\eta := -\varepsilon + \sqrt{\varepsilon^2 + 2\ell_N}$. We would then have $\kappa^2 \leq 1 + 12\ell_N$ and $\eta \leq 1/2$, at least for $n \geq 28$. Also, the exponent $-N\varepsilon\eta - \tfrac{1}{2}N\kappa^2\eta^2$ is smaller than $-\log N$, which ensures that the final, bracketed term in (11) only contributes another $O(N^{-1})$ to the $A_n(\varepsilon)$ from Theorem 1.

**5. Details of the proof for Theorem 2.** Written using the $\Psi$ function from Lemma 1, the assertion of Theorem 1 implies that

$$\Psi(z_k) = \Psi(\varepsilon\sqrt{N}) + B_n(\varepsilon) + \tau_k,$$

where, for some constant $C$,

$$B_n(\varepsilon) = N\varepsilon^4\gamma(\varepsilon) + \tfrac{1}{2}\log(1-\varepsilon^2) + \lambda_{n-k} \quad \text{and} \quad -CN^{-1} \leq \tau_k \leq C\ell_N$$

for $\varepsilon$ corresponding to the range $n/2 \geq k \leq n-1$, that is, for $0 \leq \varepsilon \leq 1 - 2N^{-1}$.

Define

$$w_k = \varepsilon\sqrt{N}\, S(\varepsilon) + \frac{\log(1-\varepsilon^2) + 2\lambda_{n-k}}{2\varepsilon\sqrt{N}\, S(\varepsilon)}.$$

We need to show that there is a constant $C'$ for which $z_k = w_k + \theta_k$, with $-C'(\varepsilon\sqrt{N}+1) \leq N\theta_k \leq C'(\varepsilon\sqrt{N}+\log N)$ for $0 \leq \varepsilon \leq 1 - 2N^{-1}$. Consider two cases.

5.1. *Suppose $\varepsilon \leq C_0/\sqrt{N}$ for some constant $C_0$.* Uniformly over that range $B_n(\varepsilon) = O(N^{-1})$ and $w_k = \varepsilon\sqrt{N} + O(N^{-1})$. From Lemma 1(i), for all nonnegative $\delta_1$ and $\delta_2$,

$$\Psi(x) + \delta_1\rho(x) \leq \Psi(x+\delta_1) \quad \text{and} \quad \Psi(x-\delta_2) + \delta_2\rho(x-\delta_2) \leq \Psi(x).$$

With $x$ equal to $\varepsilon\sqrt{N}$ and $C_1$ a large enough constant, deduce that

$$\Psi(\varepsilon\sqrt{N} - C_1N^{-1}) < \Psi(z_k) < \Psi(\varepsilon\sqrt{N} + C_1\ell_N)$$

and, hence,

$$w_k - O(N^{-1}) - C_1N^{-1} < z_k < w_k + O(N^{-1}) + C_1\ell_N.$$

5.2. *Suppose $C_0/\sqrt{N} \leq \varepsilon \leq 1 - 2N^{-1}$.* Write $x$ for $\varepsilon\sqrt{N}$ and $\beta$ for $B_n(\varepsilon) + \tau_k = \Psi(z_k) - \Psi(x)$. For all $\varepsilon$ in this range, if $C_0$ is large enough, we have $\beta > 0$ and $r(x) \leq 2/x$. The function $h(t) = t - \sqrt{t^2 + 2\beta}$ is negative, increasing and concave, with $h'(t) \leq 2\beta/t^2$. The positive numbers $\delta_1 = -h(x)$ and $\delta_2 = -h(\rho(x))$ are roots of two quadratic equations, $\delta_1 x + \tfrac{1}{2}\delta_1^2 = \beta = \delta_2\rho(x) + \tfrac{1}{2}\delta_2^2$. From Lemma 1(iii),

$$\Psi(z_k) - \Psi(x) = x\delta_1 + \tfrac{1}{2}\delta_1^2 \leq \Psi(x+\delta_1) - \Psi(x),$$

$$\Psi(x+\delta_2) - \Psi(x) \leq \rho(x)\delta_2 + \tfrac{1}{2}\delta_2^2 = \Psi(z_k) - \Psi(x),$$

which imply that $x + \delta_2 \leq z_k \leq x + \delta_1$. These bounds force $z_k$ to lie close to $x + \delta_1$:

$$0 \leq x + \delta_1 - z_k \leq \delta_1 - \delta_2 = h(\rho(x)) - h(x) \leq r(x)h'(x) \leq 4\beta/x^3 = O(\varepsilon/\sqrt{N}).$$

And $x + \delta_1$ lies close to $w_k$:

$$\begin{aligned}
x + \delta_1 &= \sqrt{N\varepsilon^2 + 2\beta} \\
&= \varepsilon\sqrt{N}\, S(\varepsilon)\left(1 + \frac{\log(1-\varepsilon^2) + 2\lambda_{n-k} + \tau_k}{N\varepsilon^2 S(\varepsilon)^2}\right)^{1/2} \\
&= w_k + \frac{\tau_k}{2\varepsilon\sqrt{N}\, S(\varepsilon)} + O(\sqrt{N}\,\ell_N^2).
\end{aligned}$$

The assertion of Theorem 2 follows.

**Acknowledgments.** We thank the referees for particularly constructive comments and suggested improvements, David Mason for explaining the role of large deviation approximations and Sandor Csörgő for providing, many years ago, an English translation of parts of Tusnády's dissertation.

Department of Statistics
and Applied Probability
University of Califonia
Santa Barbara, California 93106-3110
USA
e-mail: carter@pstat.ucsb.edu
url: www.pstat.ucsb.edu/faculty/carter/

Department of Statistics
Yale University
P.O. Box 208290
Yale Station
New Haven, Connecticut 06520-8290
USA
e-mail: david.pollard@yale.edu
url: www.stat.yale.edu/~pollard/